\documentclass[12pt,reqno]{amsart}

\usepackage{amssymb}

\usepackage{eucal}

\input{diagrams}

\font\sss=cmss8

\def\cE{{\mathcal E}}

\def\cM{{\mathcal M}}

\def\BN{{\mathbb N}}

\def\BZ{{\mathbb Z}}

\def\fm{{\mathfrak m}}

\def\fp{{\mathfrak p}}

\def\ann{\operatorname{ann}}

\def\colim{\operatorname{colim}}

\def\deg{\operatorname{deg}}
\def\depth{\operatorname{depth}}

\def\dim{\operatorname{dim}}

\def\Ext{\operatorname{Ext}}

\def\Gammam{\Gamma_{\fm}}

\def\gr{\mbox{\sf gr}\,}

\def\Gr{\mbox{\sf Gr}\,}
\def\Grsmall{\mbox{\sss Gr}}

\def\H{\operatorname{H}}
\def\hdet{\operatorname{hdet}}

\def\Hom{\operatorname{Hom}}
\def\id{\operatorname{id}}

\def\inf{\operatorname{inf}}

\def\Ker{\operatorname{Ker}}

\def\opp{\operatorname{op}}
\def\pd{\operatorname{pd}}
\def\proj{\operatorname{proj}}

\def\qgr{\mbox{\sf qgr}\,}
\def\QGr{\mbox{\sf QGr}\,}

\def\RGammam{\operatorname{R}\!\Gamma_{\fm}}

\def\RHom{\operatorname{RHom}}

\def\Spec{\operatorname{Spec}}
\def\sup{\operatorname{sup}}

\def\Tors{\mbox{\sf Tors}\,}

\numberwithin{equation}{part}

\newtheorem{Lemma}{Lemma}[section]
\newtheorem{Theorem}[Lemma]{Theorem}
\newtheorem{Proposition}[Lemma]{Proposition}

\theoremstyle{definition}
\newtheorem{Definition}[Lemma]{Definition}
\newtheorem{Setup}[Lemma]{Setup}

\newtheorem{Remark}[Lemma]{Remark}

\def\CM{Co\-hen-Ma\-cau\-lay}
\def\MCM{maximal \CM}
\def\fCMt{finite \CM\ type}

\begin{document}

\title[Non-commutative Cohen-Macaulay type]
{Finite Cohen-Macaulay type and smooth non-commutative schemes}

\author{Peter J\o rgensen}
\address{Department of Pure Mathematics, University of Leeds,
Leeds LS2 9JT, United Kingdom}
\email{popjoerg@maths.leeds.ac.uk}
\urladdr{http://www.maths.leeds.ac.uk/\~{ }popjoerg}




\begin{abstract} 

A commutative local \CM\ ring $R$ of \fCMt\ is known to be an isolated
singularity; that is, $\Spec(R) \setminus \{ \fm \}$ is smooth.

This paper proves a non-commutative analogue.  Namely, if $A$ is a
(non-commutative) graded AS \CM\ algebra which is FBN and has \fCMt,
then the non-commutative projective scheme determined by $A$ is
smooth.

\end{abstract}

\maketitle

\setcounter{section}{-1}
\section{Introduction}

Auslander proved that a commutative local \CM\ ring of \fCMt\ is an
isolated singularity.  The present paper shows a non-commutative
version of this.

To be precise about Auslander's result, consider a commutative local
noetherian ring $R$ which is Cohen-Macaulay of depth $d$.  A finitely
ge\-ne\-ra\-ted $R$-module $M$ is called \MCM\ if $\depth M = d$, and $R$ is
said to have \fCMt\ if there are only finitely many isomorphism
classes of indecomposable \MCM\ modules.  Auslander now proved that if
$R$ has \fCMt\ then it is an isolated singularity, that is, the scheme
$\Spec(R) \setminus \{ \fm \}$ is smooth.  In fact, Auslander had to
assume that $R$ was complete, see \cite[thm., p.\ 234]{Auslander}; the
general statement is due to Huneke and Leuschke, see \cite[cor.\
2]{HunekeLeuschke}.

Now, it is well known that the theory of commutative local rings has a
close analogue in the theory of non-commutative connected $\BN$-graded
algebras, so it is natural to ask for a non-commutative graded version
of Auslander's result.  Such a version is shown for FBN algebras in
theorem \ref{thm:main} below.  The notion of \fCMt\ can be carried
over directly to graded algebras, and the analogue of being an
isolated singularity is easy to guess: The canonical procedure for
removing the ``irrelevant'' maximal ideal $A_{\geq 1}$ from an
$\BN$-graded algebra $A$ is to take the projective scheme $\proj A$ of
$A$, so the analogue of $\Spec(R) \setminus \{ \fm \}$ being smooth is
that $\proj A$ is smooth.  In fact, since $A$ is non-commutative, I
must use non-commutative projective schemes of the form $\qgr A$ as
introduced in \cite{ArtinZhang}.  For $\qgr A$ to be smooth means that
each object has finite injective dimension.

The catalyst for this paper was Huneke and Leuschke's
\cite{HunekeLeuschke}, and I am inspired by their method of proof.
Thus, my basic lemma states that if $A$ has \fCMt\ and $M$ and $N$ are
\MCM\ modules, then $\Ext_A^1(M,N)$ has finite length.  However, the
proof I give of this is new:  While \cite{HunekeLeuschke} uses to good
effect that $\Ext_R^1(M,N)$ is a module over the commutative ring $R$,
there is no such aid to be had over the non-commutative ring $A$, and
so a new proof had to be found.

\begin{Setup}
\label{set:blanket}
Throughout, $k$ is a field and $A$ is a connected $\BN$-graded
noetherian $k$-algebra which is AS \CM\ in the sense that it has a
balanced dualizing complex $D$ which is equal to the $d$'th
suspension $\Sigma^d K$ of a graded $A$-bimodule $K$.
\end{Setup}

The number $d$ plays the role of dimension of $A$.  See
\cite{Yekutieli} or \cite{VdB} for information about balanced
dualizing complexes and \cite{Mori} for AS \CM\ algebras.

\medskip

Let me close the introduction with some notation which may
be convenient for the reader, although none of this differs
significantly from previous papers such as \cite{PJLocal},
\cite{PJIdent}, or \cite{JZ}.

By $\Gr A$ is denoted the category of graded $A$-left-modules and
graded homomorphisms of degree zero, and by $\gr A$ the full
subcategory of finitely generated modules.

If $M$ is in $\Gr A$ then
\[
  i(M) = \inf \{\, j \,|\, M_j \not= 0 \,\}.
\]

By $\Hom_A$ is denoted the functor
\[
  \Hom_A(-,-) = \bigoplus_{\ell} \Hom_{\Grsmall A}(-,-(\ell)),
\]
where $(\ell)$ denotes the $\ell$'th degree shift of a graded module; that
is, $(N(\ell))_j = N_{\ell+j}$.  The total right derived functor of
$\Hom_A$ is denoted $\RHom_A$, and the $i$'th derived functor of
$\Hom_A$ is denoted $\Ext_A^i$, so $\Ext_A^i \simeq \H^i\RHom_A$ and
\[
  \Ext_A^i(-,-) = \bigoplus_{\ell} \Ext_{\Grsmall A}^i(-,-(\ell)).
\]

The depth of $M$ in $\Gr A$ is defined by
\[
  \depth M = \inf \{\, i \,|\, \Ext_A^i(k,M) = 0 \,\},
\]
and $M$ in $\gr A$ is called a graded \MCM\ module if $\depth M =
d$.  Observe that $A$ itself is a graded \MCM\ module by lemma
\ref{lem:depthA} below. 

By $\Gammam$ is denoted the local section functor
\[
  \Gammam(-) = \colim \Hom_A(A/A_{\geq j},-).
\]
The total right derived functor of $\Gammam$ is denoted $\RGammam$.
By $(-)^{\prime}$ is denoted the functor $\Hom_k(-,k)$.

Finally, if $T$ is in $\Gr A$ then an element $t$ in $T$ is called
torsion if $A_{\geq j}t = 0$ for some $j$.  Graded torsion and graded
torsionfree modules are defined in the obvious way.  The full
subcategory of $\Gr A$ consisting of torsion modules is denoted by
$\Tors A$, and $\QGr A$ is defined as the quotient $\Gr A/\Tors A$,
while $\qgr A$ is defined as the full subcategory corresponding to
$\gr A$.  The quotient functor is denoted $\Gr A
\stackrel{\pi}{\longrightarrow} \QGr A$.

\section{Basic lemmas}

\begin{Lemma}
\label{lem:depthA}
The depth of $A$ is
\[
  \depth A = d.
\]
\end{Lemma}

\begin{proof}
By \cite[thm.\ 6.3]{VdB}, the algebra $A$ satisfies the so-called
$\chi$ condition of \cite[def.\ 3.7]{ArtinZhang}, so by \cite[prop.\
4.3]{PJLocal}, the depth of $M$ in $\gr A$ is
\begin{equation}
\label{equ:depth}
  \depth M
  = -\sup \{\, i \,|\, \H^i(\RGammam(M)^{\prime}) \not= 0 \,\}.
\end{equation}
In particular,
\[
  \depth A
  = -\sup \{\, i \,|\, \H^i(\RGammam(A)^{\prime}) \not= 0 \,\}
  = (*).
\]
But \cite[thm.\ 6.3]{VdB} also gives that $\RGammam(A)^{\prime}$ is
isomorphic to the balanced dualizing complex $D = \Sigma^d K$, so 
\[
  (*) = -\sup \{\, i \,|\, \H^i(\Sigma^d K) \not= 0 \,\} = d.
\]
\end{proof}

\begin{Lemma}
\label{lem:depthM}
If $M \not= 0$ is in $\gr A$ then
\[
  \depth M \leq d.
\]
\end{Lemma}

\begin{proof}
By \cite[thms.\ 6.3 and 5.1]{VdB}, I have local duality,
\[
  \RGammam(M)^{\prime} \cong \RHom_A(M,D).
\]
Along with equation \eqref{equ:depth} this gives
\[
  \depth M
  = -\sup \{\, i \,|\, \H^i(\RHom_A(M,D)) \not= 0 \,\} = (*).
\]

However, the complex $D = \Sigma^d K$ sits in cohomological
degrees $\geq -d$.  The same applies to some injective resolution of
$D$, and thus also to $\RHom_A(M,D)$.  So if $\RHom_A(M,D) \not= 0$
then
\[
  (*) \leq d
\]
proving the lemma.

To see $\RHom_A(M,D) \not= 0$, note that this follows from $M \not= 0$
because the functor $\RHom_A(-,D)$ is an equivalence of categories by
\cite[prop.\ 3.5]{Yekutieli}.
\end{proof}

\begin{Lemma}
\label{lem:extension}
Let
\[
  0 \rightarrow N \longrightarrow X \longrightarrow M \rightarrow 0
\]
be a short exact sequence in $\gr A$.  If $M$ and $N$ are graded \MCM\
modules, then so is $X$.
\end{Lemma}

\begin{proof}
There is a long exact sequence which consists of pieces
\[
  \Ext_A^i(k,N) 
  \longrightarrow \Ext_A^i(k,X) 
  \longrightarrow \Ext_A^i(k,M).
\]
This shows $\depth X \geq d$, and $\depth X \leq d$ holds by lemma
\ref{lem:depthM}.  
\end{proof}

\begin{Lemma}
\label{lem:depthadd}
Let
\[
  0 \rightarrow \Omega M \rightarrow P \rightarrow M \rightarrow 0
\]
be a short exact sequence in $\gr A$ with $P$ graded projective.  Then
\[
  0 \leq \depth M \leq d-1 \; \Rightarrow \; \depth \Omega M = \depth M + 1,
\]
and if $\Omega M \not= 0$ then
\[
  \depth M = d \; \Rightarrow \; \depth \Omega M = d.
\]
\end{Lemma}

\begin{proof}
There is a long exact sequence which consists of pieces
\[
  \Ext_A^i(k,\Omega M) \rightarrow
  \Ext_A^i(k,P) \rightarrow
  \Ext_A^i(k,M).
\]
This easily gives the first implication of the lemma because
it is clear that $\depth P = \depth A = d$.  

As for the second implication, when $\depth M = d$ the long exact
sequence implies $\depth \Omega M \geq d$, and $\depth \Omega M \leq
d$ holds by lemma \ref{lem:depthM} when $\Omega M \not= 0$.
\end{proof}

The following lemma is a special case of the graded analogue of
\cite[lem.\ 2.3]{Brown}.  I will include a proof to convince the
reader and myself, since solid references even for basic properties of
FBN rings with a grading seem hard to find.

\begin{Lemma}
\label{lem:Brown}
Let $A$ be FBN and let $M$ have minimal injective resolution $E$ in
$\Gr A$.  If $E^i$ is not torsion in the sense of the introduction,
then there exists a graded prime ideal $\fp$ of $A$ so that
\[
  \dim_k \Ext_A^i(A/\fp,M) = \infty.
\]
\end{Lemma}

\begin{proof}
By \cite[prop.\ 7.1(5)]{ArtinZhang} I have $E^i = Q \oplus T$ in $\Gr
A$ where $Q$ is a graded torsionfree injective module and $T$ a graded
torsion injective module.

For $E^i$ not to be torsion means $Q \not= 0$.  Since $E$ is minimal,
$\Ker \partial_E^i$ is graded essential in $E^i$, and so, $Q \cap \Ker
\partial_E^i$ is graded essential in $Q$.  In particular, $Q \cap \Ker
\partial_E^i \not= 0$.  It follows from \cite[lem.\ 2.1(i)]{StZ} that
there is a non-zero finitely generated graded submodule $V$ of $Q \cap
\Ker \partial_E^i$ so that $\fp = \ann_A V$ is a graded prime ideal
and so that ${}_{A/\fp}V$ is non-singular (in the ungraded sense).
Let $e$ be a non-zero graded element of $V$ so
\[
  0 \not= e \in V \subseteq Q \cap \Ker \partial_E^i \subseteq E^i.
\]

Now, $V \subseteq Q$ implies that $V$ is graded torsionfree.  Since
$V$ is annihilated by $\fp$ but is non-zero and torsionfree, $A_{\geq
1}$ cannot be contained in $\fp$.  So $(A/\fp)_{\geq 1}$ is a non-zero
ideal of $A/\fp$, and hence there is a regular graded element $c$ of
positive degree in $A/\fp$; this follows e.g.\ from \cite[lem.\
2.1(iii)]{StZ}.

Now note that $X$, the subcomplex of $E$ consisting of elements annihilated
by $\fp$, is isomorphic to $\Hom_A(A/\fp,E)$, and that hence,
\[
  \H^i\!X \cong \H^i \Hom_A(A/\fp,E) \cong \Ext_A^i(A/\fp,M).
\]
Consider $c^m e$ for $m \geq 0$; these are elements of $E^i$.  In
fact, they are elements of $V$, and this means that they are
annihilated by $\fp$ so are elements of $X^i$.  It also means that
they are in $\Ker \partial_E^i$, hence in $\Ker \partial_X^i$, so
represent classes in $H^i X \cong \Ext_A^i(A/\fp,M)$.  

If these classes are non-zero then they must be different because the
$c^m e$ have different graded degrees, and so $\dim_k
\Ext_A^i(A/\fp,M) = \infty$ as desired.  So to finish the proof, I
must see for each $m$ that $c^m e$ does not represent zero in $H^i X$.
Suppose to the contrary that it does for some $m$.  Then there is an
$x$ in $X^{i-1}$ with $c^m e = \partial_X^{i-1}(x)$.  Hence
\begin{align}
\nonumber
  (A/\fp) c^m e 
    & = (A/\fp) \partial_X^{i-1}(x)
    = \partial_X^{i-1}((A/\fp)x)
    \cong \frac{(A/\fp)x}{(A/\fp)x \cap \Ker \partial_X^{i-1}} \\
\label{equ:f}
    & = \frac{(A/\fp)x}{(A/\fp)x \cap \Ker \partial_E^{i-1}},
\end{align}
where the last $=$ is because $X$ is just a subcomplex of $E$.  But
$\Ker \partial_E^{i-1}$ is graded essential in $E^{i-1}$ so $(A/\fp)x
\cap \Ker \partial_E^{i-1}$ is graded essential in $(A/\fp)x$, and by
\cite[lem.\ A.I.2.8]{NV} this implies that $(A/\fp)x \cap \Ker
\partial_E^{i-1}$ is essential (in the ungraded sense) in $(A/\fp)x$.
Hence the last quotient module in equation \eqref{equ:f} is singular
over $A/\fp$ by \cite[prop.\ 3.26]{GW} so the equation shows that
$(A/\fp)c^m e$ is singular over $A/\fp$.

However, this module is contained in the non-singular module
${}_{A/\fp}V$ so must be zero, whence $c^m e = 0$.  As $c^m$ is
regular in $A/\fp$ this implies by \cite[prop.\ 6.9]{GW} that $e$ is
in the singular submodule of ${}_{A/\fp}V$.  But $e$ is non-zero so
this shows that ${}_{A/\fp}V$ cannot be non-singular; a contradiction.
\end{proof}

\section{Finite \CM\ type and smoothness}

\begin{Remark}
\label{rmk:KrullSchmidt}
The category $\gr A$ is a $k$-linear category with finite dimensional
$\Hom$ spaces.  This implies that $\gr A$ is a Krull-Schmidt
ca\-te\-go\-ry.  That is, each object is a direct sum of finitely many
uniquely determined indecomposable objects.

Note that if $M$ decomposes as $M \langle 1 \rangle \oplus \cdots
\oplus M \langle s \rangle$ in $\gr A$ then $M$ is graded \MCM\ if and
only if each $M \langle m \rangle$ is graded \MCM.
\end{Remark}

\begin{Definition}
\label{def:fCMt}
The algebra $A$ is said to have \fCMt\ if there exist finitely many
indecomposable graded \MCM\ modules $Z \langle 1 \rangle, \ldots, Z
\langle t \rangle$ so that, up to isomorphism, the indecomposable
graded \MCM\ modules in $\gr A$ are precisely the degree shifts $Z
\langle n \rangle(\ell)$ for $1 \leq n \leq t$ and $\ell \in \BZ$.
\end{Definition}

Each $Z \langle n \rangle$ can clearly be replaced with any degree
shift, and so if convenient, I can suppose 
\[
  i(Z \langle n \rangle) = 0
\]
for each $n$.

\begin{Lemma}
\label{lem:Auslander}
Let $A$ have \fCMt\ and let $M$ and $N$ in $\gr A$ be graded \MCM\
modules.  Then
\[
  \dim_k \Ext_A^1(M,N) < \infty.
\]
\end{Lemma}

\begin{proof}
Without loss of generality, I can suppose that $M$ is indecomposable
and that $N$ sits in graded degrees $\geq 0$.

Using a free resolution of $M$ in $\gr A$ easily shows
\[
  \Ext_{\Grsmall A}^1(M,N(\ell)) = 0 \; \mbox{for} \; \ell \ll 0.
\]
Using also that $\gr A$ has finite dimensional $\Hom$ spaces
shows
\[
  \dim_k \Ext_{\Grsmall A}^1(M,N(\ell)) < \infty
\]
for each $\ell$.  Since
\[
  \Ext_A^1(M,N) = \bigoplus_{\ell} \Ext_{\Grsmall A}^1(M,N(\ell)),
\]
the lemma will follow if I can show
\[
  \Ext_{\Grsmall A}^1(M,N(\ell)) = 0 \; \mbox{for} \; \ell \gg 0.
\]
That is, I must show that for $\ell \gg 0$, any short exact sequence
\begin{equation}
\label{equ:c}
  0 \rightarrow N(\ell) \longrightarrow X \longrightarrow M \rightarrow 0
\end{equation}
in $\gr A$ is split.  

Observe that in such a sequence, $X$ is graded \MCM\ by lemma
\ref{lem:extension}.  Hence 
\[
  X \cong \bigoplus_m X \langle m \rangle
\]
where the $X \langle m \rangle$ are indecomposable graded \MCM\
modules.  Since $N$ sits in graded degrees $\geq 0$, it can be
generated by graded elements of degrees $0$ to $g$ for some $g \geq
0$.  Let
\[
  X^{\prime} 
    = \bigoplus_{i(X \langle m \rangle) \leq g-\ell} X \langle m \rangle 
  \;\;\mbox{and}\;\;
  X^{\prime\prime} 
    = \bigoplus_{i(X \langle m \rangle) > g-\ell} X \langle m \rangle
\]
so that $i(X^{\prime\prime}) > g - \ell$ and
\[
  X \cong X^{\prime} \oplus X^{\prime\prime}.
\]
Note that $X$ and the $X \langle m \rangle$ depend on $\ell$.  Hence
I cannot conclude $X^{\prime} = 0$ for $\ell \gg 0$; indeed, this
turns out to be false.

The homomorphism $N(\ell) \longrightarrow X$ in \eqref{equ:c} consists
of components 
\[
  N(\ell) \longrightarrow X^{\prime}
  \;\; \mbox{and} \;\;
  N(\ell) \longrightarrow X^{\prime\prime}.  
\]
Since $N$ can be generated by graded elements of degrees $0$ to $g$,
it follows that $N(\ell)$ can be generated by elements of degrees
$-\ell$ to $g-\ell$, so the homomorphism $N(\ell) \longrightarrow
X^{\prime\prime}$ is zero because $i(X^{\prime\prime}) > g - \ell$.
So $N(\ell) \longrightarrow X$ factors through the
inclusion $X^{\prime} \hookrightarrow X$, and this means that there is
a commutative diagram with exact rows,
\[
  \begin{diagram}[labelstyle=\scriptstyle,width=1cm,height=1cm]
    0 & \rTo & N(\ell)
      & \rTo & X^{\prime}
      & \rTo & M^{\prime}
      & \rTo & 0 \\
      &      & \vEq
      &      & \dTo
      &      & \dTo
      &      & \\
    0 & \rTo & N(\ell)
      & \rTo & X
      & \rTo & M
      & \rTo & 0 \lefteqn{.} \\
  \end{diagram}
\]
Applying the Snake Lemma embeds this into a diagram with exact
rows and colums,
\begin{equation}
\label{equ:b}
  \begin{diagram}[labelstyle=\scriptstyle,width=1cm,height=1cm]
      &      &
      &      & 0
      &      & 0
      &      & \\
      &      &
      &      & \dTo
      &      & \dTo
      &      & \\
    0 & \rTo & N(\ell)
      & \rTo & X^{\prime}
      & \rTo & M^{\prime}
      & \rTo & 0 \\
      &      & \vEq
      &      & \dTo
      &      & \dTo
      &      & \\
    0 & \rTo & N(\ell)
      & \rTo & X
      & \rTo & M
      & \rTo & 0 \lefteqn{.} \\
      &      &
      &      & \dTo
      &      & \dTo
      &      & \\
      &      &
      &      & X^{\prime\prime}
      & \rEq & X^{\prime\prime}
      &      & \\
      &      &
      &      & \dTo
      &      & \dTo
      &      & \\
      &      &
      &      & 0
      &      & 0
      &      & \\
  \end{diagram}
\end{equation}
The first vertical exact sequence is split by construction, and this
implies that the second vertical exact sequence is split.  But
$M$ is indecomposable, so either $M^{\prime} = 0$, $X^{\prime\prime}
\cong M$ or $M^{\prime} \cong M$, $X^{\prime\prime} = 0$.  

Let me now assume $\ell \gg 0$ and show that
$M^{\prime} = 0$, $X^{\prime\prime} \cong M$.  Suppose to the contrary
that $X^{\prime\prime} = 0$.  Then $X = X^{\prime}$ so
\begin{equation}
\label{equ:i}
  X
  = \bigoplus_{i(X \langle m \rangle) \leq g-\ell} X \langle m \rangle.
\end{equation}

Each $X \langle m \rangle$ is an indecomposable graded \MCM\
module, so each $X \langle m \rangle$ is a degree shift of one of
the $Z\langle n \rangle$ from definition \ref{def:fCMt},
\begin{equation}
\label{equ:a}
  X \langle m \rangle \cong Z \langle n \rangle (p).
\end{equation}
As remarked after that definition, I can suppose $i(Z \langle n
\rangle) = 0$ for each $n$, and then $i(X \langle m \rangle) \leq
g-\ell$ implies $p \geq \ell-g$.  Since there are only finitely many
$Z \langle n \rangle$'s, there is a $q$ such that each $Z \langle n
\rangle$ is generated by graded elements of degrees $\leq q$,
and then $p \geq \ell-g$ means that each $Z \langle n \rangle (p)$ is
generated by elements of degrees $\leq q+g-\ell$.  By equation
\eqref{equ:a}, the same holds for each $X \langle m \rangle$, and
since $\ell \gg 0$ implies $q+g-\ell < i(M)$, each homomorphism $X
\langle m \rangle \longrightarrow M$ must be zero.

But there is a surjection $X \longrightarrow M$ which by equation
\eqref{equ:i} must then also be zero, a contradiction.  Hence
$X^{\prime\prime} \not= 0$, and $M^{\prime} = 0$, $X^{\prime\prime}
\cong M$ holds as claimed.

The first horizontal exact sequence in diagram
\eqref{equ:b} then shows $X^{\prime} \cong N(\ell)$, and so the
original exact sequence \eqref{equ:c} reads
\begin{equation}
\label{equ:d}
  0 \rightarrow N(\ell) \longrightarrow N(\ell) \oplus M
  \longrightarrow M \rightarrow 0.
\end{equation}
I have not yet proved that the sequence is split, since I have not
identified the homomorphisms.  However, $N(\ell) \longrightarrow
N(\ell) \oplus M$ consists of components 
\[
  N(\ell) \longrightarrow N(\ell)
  \;\; \mbox{and} \,\, 
  N(\ell) \longrightarrow M, 
\]
and since $\ell \gg 0$ implies $g - \ell < i(M)$, the homomorphism
$N(\ell) \longrightarrow M$ is zero because $N$ can be generated by
graded elements of degrees $0$ to $g$ and $N(\ell)$ by elements of
degrees $-\ell$ to $g-\ell$.  Hence $N(\ell) \longrightarrow N(\ell)$
must be injective, and since $N(\ell)$ is finitely generated, each of
its graded components is finite dimensional over $k$ by \cite[prop.\
2.1]{ArtinZhang}, so it follows that $N(\ell) \longrightarrow N(\ell)$
is also surjective.  Hence $N(\ell) \longrightarrow N(\ell)$ is
bijective so there is a splitting of $N(\ell) \longrightarrow N(\ell)
\oplus M$, and this proves that \eqref{equ:d} and hence \eqref{equ:c}
is split as desired.  
\end{proof}

\begin{Proposition}
\label{pro:injectiveresolution}
Let $A$ be FBN with \fCMt\ and let $M$ in $\gr A$ have minimal
injective resolution $E$ in $\Gr A$.  Then $E^d, E^{d+1}, \ldots$ are
torsion.
\end{Proposition}

\begin{proof}
This is clear for $M = 0$.  Let me next give a proof when $M$ is
graded \MCM.  Suppose that $E^{d+i}$ is non-torsion for some $i \geq
0$.  By lemma \ref{lem:Brown}, there exists a graded prime ideal $\fp$
with 
\[
  \dim_k \Ext_A^{d+i}(A/\fp,M) = \infty.
\]

This implies that $\fp$ is not the maximal ideal $A_{\geq 1}$, since,
as noted earlier, $A$ satisfies condition $\chi$ of \cite[def.\
3.7]{ArtinZhang} by \cite[thm.\ 6.3]{VdB}.  Hence $\depth A/\fp \geq
1$, for if $\depth A/\fp = 0$ then there would exist a non-zero
homomorphism $k(\ell) \stackrel{\varphi}{\longrightarrow} A/\fp$, and
this would lead to a contradiction.  Namely, $\varphi(1) = x \not= 0$
would imply that $x$ was represented by $y \in A \setminus \fp$, and
$A_{\geq 1} x = A_{\geq 1} \varphi(1) = \varphi(A_{\geq 1} 1) =
\varphi(0) = 0$ would imply $A_{\geq 1} y \subseteq \fp$ and hence
$A_{\geq 1} \cdot Ay \subseteq \fp$.  As $\fp$ is a prime ideal, this
would mean either $A_{\geq 1} \subseteq \fp$, contradicting that $\fp$
is not $A_{\geq 1}$, or $Ay \subseteq \fp$, contradicting $y \in A
\setminus \fp$. 

Note that $\depth A/\fp \geq 1$ implies $d \geq 1$ by lemma 
\ref{lem:depthM}.  Now let
\[
  0 \rightarrow \Omega^{d+i-1}(A/\fp)
    \longrightarrow P_{d+i-2}
    \longrightarrow \cdots
    \longrightarrow P_0 
    \rightarrow A/\fp \rightarrow 0
\]
be an exact sequence in $\gr A$ where the $P_j$ are graded projective.
This clearly gives $\Ext_A^{d+i}(A/\fp,M) \cong
\Ext_A^1(\Omega^{d+i-1}(A/\fp),M)$ so
\[
  \dim_k \Ext_A^1(\Omega^{d+i-1}(A/\fp),M) = \infty.
\]
Hence $\Omega^{d+i-1}(A/\fp)$ cannot be zero, and since $\depth A/\fp
\geq 1$, lemma \ref{lem:depthadd} implies that $\Omega^{d+i-1}(A/\fp)$
is graded \MCM.  Lemma \ref{lem:Auslander} thus says
\[
  \dim_k \Ext_A^1(\Omega^{d+i-1}(A/\fp),M) < \infty,
\]
contradiction the previous equation.  So $E^d, E^{d+1}, \ldots$ are
torsion. 

Now let $M$ be any finitely generated graded module.  Let
\begin{equation}
\label{equ:e}
  0 \rightarrow \Omega^d M
    \longrightarrow Q_{d-1}
    \longrightarrow \cdots
    \longrightarrow Q_0 
    \rightarrow M \rightarrow 0
\end{equation}
be an exact sequence in $\gr A$ where the $Q_j$ are graded projective.  Lemma
\ref{lem:depthadd} implies that $\Omega^d M$ is either $0$ or graded
\MCM, and the $Q_j$ are also either $0$ or graded \MCM, so I have
already proved that the proposition applies to all these modules.  So
I can prove the proposition for $M$ by working along the sequence
\eqref{equ:e} from the left hand end, using the following fact: If
\[
  0 \rightarrow K
    \longrightarrow Q
    \longrightarrow L
    \rightarrow 0
\]
is a short exact sequence in $\gr A$ where the proposition applies to
$K$ and $Q$, then it also applies to $L$.

To prove this, let $E_K$ and $E_Q$ be the minimal injective
resolutions of $K$ and $Q$ in $\Gr A$.  The homomorphism $K
\longrightarrow Q$ induces a chain map $E_K \longrightarrow E_Q$; let
$E$ be the mapping cone.  The long exact cohomology sequence shows
$\H(E) \cong L$ and by construction $E$ has the form
\[
  E = \cdots \rightarrow 0 
             \rightarrow E_K^0
             \rightarrow E_K^1 \oplus E_Q^0
             \rightarrow E_K^2 \oplus E_Q^1
             \rightarrow \cdots;
\]
here $E_K^0$ is in cohomological degree $-1$ and splits away, so all
in all there is an injective resolution of $L$ of the form
\[
  \widetilde{E}
  = \cdots \rightarrow 0 
           \rightarrow (E_K^1 \oplus E_Q^0)/E_K^0
           \rightarrow E_K^2 \oplus E_Q^1
           \rightarrow E_K^3 \oplus E_Q^2
           \rightarrow \cdots.
\]
If the proposition applies to $K$ and $Q$, then the modules $E_K^d,
E_K^{d+1}, \ldots$ and $E_Q^d, E_Q^{d+1}, \ldots$ are torsion.  Then
$\widetilde{E}^d, \widetilde{E}^{d+1}, \ldots$ are also torsion, and
since the minimal injective resolution $E_L$ of $L$ in $\Gr A$ is a
direct summand in any injective resolution of $L$ in $\Gr A$, and so a
direct summand in $\widetilde{E}$, it follows that $E_L^d, E_L^{d+1},
\ldots$ are torsion.
\end{proof}

The following is the main result of this paper.

\begin{Theorem}
\label{thm:main}
Recall the standing setup \ref{set:blanket}.  Let $A$ be FBN with
\fCMt.  Then each $\cM$ in $\qgr A$ has $\id \cM \leq d - 1$.
\end{Theorem}

\begin{proof} 
Consider the quotient functor $\Gr A \stackrel{\pi}{\longrightarrow}
\QGr A$.  There exists an $M$ in $\gr A$ with $\cM \cong \pi(M)$ (see
\cite[p.\ 234]{ArtinZhang}), and if $E$ is a minimal injective
resolution of $M$ in $\Gr A$ then it follows easily from \cite[prop.\
7.1]{ArtinZhang} that $\cE = \pi(E)$ is an injective resolution of
$\cM$.  But $E^d, E^{d+1}, \ldots$ are torsion by proposition
\ref{pro:injectiveresolution}, so $\cE^d = \cE^{d+1} = \cdots = 0$.
\end{proof}

\section{An example}

This section contains an example inspired by
\cite{AuslanderRational}. 

Suppose that the ground field $k$ does not have characteristic $2$,
but contains a primitive $n$'th root of unity, $q$, and let $B =
k\langle x,y \rangle / (yx - qxy)$ where $x$ and $y$ have degree $1$.
Let $G = \langle g \rangle$ be the cyclic group of order $2$, and let
$G$ act on $B$ by $gb = (-1)^{\deg b}b$.  

It is clear that the fixed ring $A = B^G$ is given by
\[
  A_j
  = \left\{
      \begin{array}{cl}
        B_j & \mbox{for $j$ even,} \\
        0   & \mbox{for $j$ odd.} \\
      \end{array}
    \right.
\]
I will show that $A$...
\begin{enumerate}

  \item  falls under setup \ref{set:blanket} with $d=2$, 

\smallskip

  \item  has \fCMt, 

\smallskip

  \item  is FBN,

\smallskip

  \item  has infinite global dimension.

\end{enumerate}
These properties imply that theorem \ref{thm:main} applies to $A$
and says that each $\cM$ in $\qgr A$ has $\id \cM \leq 1$, although
$A$ itself has infinite global dimension.

To prove (i) through (iv), let me first consider $B$ more closely.  It
is a so-called twist of the commutative polynomial ring $k[x,y]$ by
the algebraic twisting system given by
\[
  \tau_n(x) = x, \;\; \tau_n(y) = q^{-n}y, 
\]
cf.\ \cite[exam.\ 3.6]{Zhang}, and so the results of \cite{Zhang}
imply that it is a connected $\BN$-graded noetherian $k$-algebra which
is AS regular of global dimension $2$.  

By \cite[sec.\ 3]{JZ}, the fixed ring $A = B^G$ is therefore a
connected $\BN$-graded noetherian $k$-algebra with a balanced
dualizing complex, and \cite[lem.\ 3.1]{JZ} shows that the balanced
dualizing complex is concentrated in cohomological degree $-2$,
proving (i).

In fact, it is not hard to check that if I denote by $\hdet$ the
homological determinant defined in \cite[sec.\ 2]{JZ}, then $\hdet(g)
= 1$ whence $A = B^G$ is AS Gorenstein by \cite[thm.\ 3.3]{JZ}.  This
implies that $A$ is a dualizing complex for itself, and hence, since
$d = 2$, it follows from \cite[thm.\ 3.9 and cor.\ 4.10]{Yekutieli}
that there is an automorphism $\sigma$ so that the second suspension
$\Sigma^2(A^{\sigma})$ is a balanced dualizing complex.  This will be
handy below for the proof of (ii).

It is a consequence of \cite[prop.\ 5.6(a)]{Zhang} that $B$ is a PI
algebra, so the subalgebra $A$ is also PI and hence FBN, proving (iii).

The Hilbert series of $A$ is
\[
  H_A(t) = 1 + 3t^2 + 5t^4 + \cdots = \frac{1+t^2}{(1-t^2)^2},
\]
and as this is not $1$ divided by a polynomial, $A$ cannot have finite
global dimension, proving (iv).

Finally, to prove (ii), let $M$ in $\gr A$ be an indecomposable graded
\MCM\ module.  Since $\Sigma^2(A^{\sigma})$ is a balanced dualizing
complex, $A^{\sigma}$ is a balanced dualizing module in the
terminology of \cite[sec.\ 4]{Mori}.  By \cite[lem.\ 4.6]{Mori} I
therefore have 
\[
  M \cong \Hom_{A^{\opp}}(\Hom_A(M,A^{\sigma}),A^{\sigma}), 
\]
and the $\sigma$'s cancel so
\[
  M \cong \Hom_{A^{\opp}}(\Hom_A(M,A),A).
\]

Writing $N = \Hom_A(M,A)$ thus gives $M \cong \Hom_{A^{\opp}}(N,A)$.
The inclusion $A \hookrightarrow B$ is split when viewed as a
homomorphism over $A^{\opp}$, cf.\ \cite[sec.\ 3]{JZ}, and so there is
a split inclusion over $A$,
\begin{equation}
\label{equ:h}
  M \cong \Hom_{A^{\opp}}(N,A) \hookrightarrow \Hom_{A^{\opp}}(N,B). 
\end{equation}

Let me show that $\Hom_{A^{\opp}}(N,B)$ is in fact a graded projective
$B$-left-module.  It is the zeroth cohomology of
$\RHom_{A^{\opp}}(N,B)$, so there is a distinguished triangle over
$B$,
\[
  \Hom_{A^{\opp}}(N,B) 
  \longrightarrow \RHom_{A^{\opp}}(N,B)
  \longrightarrow X
  \longrightarrow,
\]
where the cohomology of $X$ is concentrated in cohomological degrees
$\geq 1$.  This again gives a distinguished triangle
\begin{align}
\label{equ:g}
\nonumber
  \RHom_B(k,\Hom_{A^{\opp}}(N,B))
  & \longrightarrow \RHom_B(k,\RHom_{A^{\opp}}(N,B)) \\
  & \;\;\;\;\;\;\;\;\;\;\;\;\;\;\;\;\;\;\; \longrightarrow \RHom_B(k,X) \longrightarrow.
\end{align}
The minimal free resolution of $k$ over $B$ from \cite[exam.\
3.6]{Zhang} implies
\[
  \Ext_B^i(k,B)
  = \left\{
      \begin{array}{cl}
        k(2) & \mbox{for $i=2$,} \\
        0    & \mbox{otherwise;} \\
      \end{array}
    \right.
\]
in other words $\RHom_B(k,B) = \Sigma^{-2}k(2)$, and so the second
term of \eqref{equ:g} is
\begin{align*}
  \RHom_B(k,\RHom_{A^{\opp}}(N,B)) 
  & \cong \RHom_{A^{\opp}}(N,\RHom_B(k,B)) \\
  & \cong \RHom_{A^{\opp}}(N,\Sigma^{-2}k(2)) \\
  & \cong \Sigma^{-2}\RHom_{A^{\opp}}(N,k)(2);
\end{align*}
the cohomology of this is concentrated in cohomological degrees $\geq
2$.  Moreover, since the cohomology of $X$ is concentrated in
cohomological degrees $\geq 1$, the same applies to the cohomology of
the third term of \eqref{equ:g},
\[
  \RHom_B(k,X).
\]

The cohomology long exact sequence of \eqref{equ:g} now shows that the
cohomology of the first term of the triangle,
\[
  \RHom_B(k,\Hom_{A^{\opp}}(N,B)),
\]
is concentrated in cohomological degrees $\geq 2$, that is, 
\[
  \depth_B \Hom_{A^{\opp}}(N,B) \geq 2.  
\]
But then the Auslander-Buchsbaum formula, \cite[thm.\ 3.2]{PJIdent},
implies
\[
  \pd_B \Hom_{A^{\opp}}(N,B) \leq 0, 
\]
and so $\Hom_{A^{\opp}}(N,B)$ is graded projective.

So the split inclusion \eqref{equ:h} says that the graded
$A$-left-module $M$ is a direct summand of some graded projective
$B$-left-module, viewed as an $A$-left-module.  This again implies
that $M$ is a direct summand of some graded free $B$-left-module,
$B(\ell_1) \oplus \cdots \oplus B(\ell_s)$, viewed as an
$A$-left-module.  Given that $M$ is indecomposable, it is already a
direct summand in one of the $B(\ell_j)$, and so, if I decompose $B$
as an $A$-left-module, the resulting direct summands are, up to degree
shift, the only possible indecomposable graded \MCM\ modules in $\gr
A$.

As $B$ is finitely generated over $A$ by \cite{JZ}, there are only
finitely many direct summands, so $A$ has \fCMt\ proving (ii).

\bigskip
\noindent
{\bf Acknowledgement. }  The diagrams were typeset with Paul Taylor's
{\tt diagrams.tex}.

\end{document}